\documentclass[a4paper,10pt]{article}
\usepackage{amssymb,amsmath}
\newcommand{\Q}{\mathbb{Q}}
\newcommand{\F}{\mathbb{F}}

\renewcommand{\b}[1]{{\bf #1}}
\newcommand{\x}{\b{x}}
\newcommand{\be}{\b{e}}
\newcommand{\cl}[1]{{\cal #1}}
\newcommand{\Z}{\mathbb{Z}}
\newcommand{\lb}[1]{\mbox{\scriptsize\bf #1}}

\newcommand{\ep}{\varepsilon}

\newtheorem{theorem}{Theorem}
\newtheorem{lemma}{Lemma}
\newtheorem{question}{Question}

\begin{document}
\title{Zeros of $p$-adic forms}
\author{D.R. Heath-Brown\\Mathematical Institute, Oxford}
\date{}
\maketitle
\section{Introduction}

This paper will be concerned with the existence of $p$-adic zeros of
$p$-adic forms.  We shall be concerned mainly, but not solely, with
quartic forms.  Before stating our results it is appropriate to recall
the general situation.   

Let $p$ be a prime and let
$F(x_1,\ldots,x_n)=F(\x)\in\Q_p[x_1,\ldots,x_n]$ be a form of degree
$d$.  It was
conjectured by Artin \cite[p.x]{Artin} that $\Q_p$ is a $C_2$ field, so that
$F(\x)$ should have a non-trivial $p$-adic zero as soon as $n>d^2$.
It is fairly easy to construct examples for every $p$ and every $d$ in
which $n=d^2$ and $F(\x)$ has no non-trivial $p$-adic zero.  It
follows readily from work of 
Brauer \cite{Br} that for every $d$ there is a number $v_d$
such that, for every $p$, the form $F(\x)$ has a non-trivial $p$-adic
zero as soon as $n>v_d$.  Brauer's method was elementary, and used
multiply nested inductions.  The resulting value of $v_d$ was too
large to write down.  Ax and Kochen \cite{AK} used methods from
mathematical logic to show that for every $d$ there is a number $p(d)$
such that every form with $n>d^2$ and $p>p(d)$ has a non-trivial
$p$-adic zero.  Later work by  Brown \cite{Brown} provided a value for
$p(d)$.  If one writes $a\uparrow b$ for $a^b$ then he showed one could take  
\begin{equation}\label{bb}
p(d)=2\uparrow(2\uparrow(2\uparrow(2\uparrow(2\uparrow(d\uparrow(11\uparrow
4d)))))).
\end{equation}

In the opposite direction, Terjanian \cite{Ter} showed that Artin's
conjecture is false in general, by providing a counterexample with
$p=2, d=4$ and $n=18$.  Later work, by Lewis and Montgomery \cite{LM}
amongst others, gives many more counterexamples.  In particular
\cite[Theorem 1]{LM} shows that for every $p$ and every  $\ep>0$ there
are infinitely many degrees $d$ and corresponding forms $F$
with no $p$-adic zero, and for which
\[n>\exp\{\frac{d}{(\log d)(\log\log d)^{1+\ep}}\}.\]
It should be noted however that all the known counterexamples to
Artin's conjecture have even degree $d$.

Since the original conjecture of Artin is now known to be false, the
natural questions become:- 
\begin{enumerate}
\item[(1)] For which values of $d$ is Artin's
conjecture true?
\item[(2)] How small can one take $v_d$ in Brauer's theorem? 
\item[(3)] How small can one
take $p(d)$ in the Ax-Kochen theorem?  
\item[(4)] What can one say about values
$v_d(p)$ for which every $p$-adic form of degree $d$ in $n>v_d(p)$
variables has a non-trivial zero?
\end{enumerate}
As to the first question, it is classical that Artin's conjecture
holds for degree $2$, and the case $d=3$ was handled by Lewis
\cite{Lewis}.  Thus the first case of interest is that of degree 4.

Turning to the number $p(d)$ in the Ax-Kochen theorem, another result
of Ax and Kochen \cite{AK2}
shows that the theory of $p$-adic fields is decidable.
Thus for each fixed prime $p$ and each fixed degree $d$ there is, in
principle, a procedure for deciding the truth or otherwise of the
statement:-

Every form $F(x_1,\ldots,x_{d^2+1})\in\Q_p[x_1,\ldots,x_{d^2+1}]$ has
a nontrivial zero over $\Q_p$. 

It follows that one can, in theory, test each prime up to Brown's
bound (\ref{bb}), 
thereby deciding whether or not Artin's conjecture holds for a 
particular degree 
$d$.  A more practical approach has its origins in the work of  Lewis
\cite{Lewis} (for $d=3$), of Birch and Lewis \cite{BL} (for $d=5$), 
and of Laxton 
and Lewis \cite{LL} (for $d=7$ and 11).  These papers consider forms
over $\Z_p$ and their reductions 
modulo $p$.  Provided that $n>d^2$, a $p$-adic reduction argument followed 
by an application of the Chevalley--Warning Theorem
produces a form modulo $p$ with a non-trivial zero.  By Hensel's
Lemma, if this zero 
is non-singular modulo $p$ it can be lifted to a 
$p$-adic zero.  Thus the crux of the problem is to find non-singular
zeros modulo $p$.  Lewis's 
argument resolved this for all $p$ when $d=3$, but in the other cases
the method only works for 
sufficiently large $p$.  Moreover the method appears to break down 
completely if the degree $d$ 
is composite, or can be written as a sum of composite numbers.  
Thus if $d=4$, for example, and $p$ is an odd prime for which $\nu$,
say, is a quadratic non-residue, one can construct forms
\begin{equation}\label{bad}
(x_1^2+\ldots+x_{n-1}^2)^2-\nu x_n^4
\end{equation}
in an arbitrary number of variables, but which have no non-singular
zero modulo $p$.

However, in those cases
where the method is successful, it can be adapted to provide
reasonable values for $p(d)$.  In particular Leep and
Yeomans \cite{LY}
show that if $d=5$ then Artin's 
conjecture holds for all primes $p\ge 47$.  Our first result gives a small
improvement on this.
\begin{theorem}\label{Th5}
Let $F(x_1,\ldots,x_n)=F(\x)\in\Q_p[x_1,\ldots,x_n]$ be a form of degree $5$
with $n>25$.  Then if
$p\ge 17$ there is a non-zero vector $\x\in\Q_p^n$ with $F(\x)=0$.
\end{theorem}

While our method fails for $p\le 13$ there is a variant of it which
might work at least for some such primes.  Since moreover we know of
no counterexamples to Artin's conjecture for $d=5$, we ask the
following question.
\begin{question}
Does Artin's conjecture hold for $d=5$ for all primes $p$?
\end{question}

In situations where the above approach fails, and in particular for
quartic forms, we can only handle small primes by versions of Brauer's
argument. The basic idea is to show via an induction argument that
$F(\x)$ represents a diagonal form in a reasonably large number of
variables.  Thus one finds linearly independent vectors
$\be_1,\ldots,\be_m\in\Q_p^n$ such that
\[F(t_1\be_1+\ldots +t_m\be_m)=\sum_{i=1}^m c_m t_m^d.\]
In general this will only be possible when $n$ is very much larger
than $m$.  However existence questions for $p$-adic zeros of diagonal
forms are relatively routine, and one can show that a non-trivial zero
always exists when $m>d^2$ (see Davenport and Lewis \cite{DL}), 
and often for smaller $m$.

Quasi-diagonalization techniques have been refined by various 
authors, and work of Wooley
\cite{W} gives the best general bounds currently available.  In
particular Wooley \cite[Corollary 1.1]{W} shows that we may take
\begin{equation}\label{fgh}
v_d\le d^{2^d}.
\end{equation}
(Recall that every $p$-adic form of degree $d$, in $n$ variables,
has a non-trivial $p$-adic
zero as soon as $n>v_d$.)  In particular we have
\begin{equation}\label{32}
v_4\le 2^{32}=4294967296.
\end{equation}
However (\ref{fgh}) is intended merely as a neat expression, valid for
all $d$, and Wooley's analysis gives more accurate information if we
specialize to $d=4$, as we shall describe in \S \ref{wooleysec}.  This
leads to the bound
\begin{equation}\label{med}
v_4\le 623426.
\end{equation}

This is a considerable improvement on (\ref{32}), and it is
in the context of this better estimate that our new bounds should be
judged. We shall prove the following results.
\begin{theorem}\label{odd}
We have 
\begin{enumerate}
\item[(i)] $v_4(p)\le 128$ for $p=3$ and $p=7$;
\item[(ii)] $v_4(5)\le 312$;
\item[(iii)] $v_4(p)\le 120$ for $p\ge 11$.
\end{enumerate}
\end{theorem}
\begin{theorem}\label{even}
We have $v_4(2)\le 9126$.
\end{theorem}
The case $p=2$ has been stated separately since it transpires that a rather
different approach is required in this case.

For Theorem \ref{odd} the technique we shall adopt is a hybrid between
Brauer's quasi-diagonalization procedure and the
$p$-adic reduction method.  In particular we shall not reduce $F(\x)$
to a completely diagonal shape, but instead produce a form whose
reduction modulo $p$ can be guaranteed to have a non-singular zero.
One cannot do this without forcing certain coefficients to vanish, as
examples of the shape (\ref{bad}) demonstrate.  Thus instead of
producing a form which is diagonal, we merely produce one whose
reduction modulo $p$ avoids certain excluded types.

Our analysis of Theorem \ref{odd}, and also to a lesser extent that of
Theorem \ref{even}, can be viewed as reducing the problem to one in
which we have to solve a system of $k$ simultaneous quadratic forms in
$m$ variables over $\Q_p$.  We write $\beta(k;\Q_p)$ for the largest
$m$ for which there is such a system with no non-trivial common zero
over $\Q_p$.  Then Lemma \ref{Ll} shows that $v_4(p)\le
\beta(8;\Q_p)+16$ for $p\not=2,5$, for example, while Lemma \ref{lb} shows that
$v_4(p)\ge\beta(4;\Q_p)$. Thus it is natural to ask what one would
expect to be the true size of $\beta(r;\Q_p)$.  Artin's original
conjecture implies that $\beta(r;\Q_p)=4r$ for all $p$, and the
Ax-Kochen theorem shows that this holds for $p\ge p(r)$. It is classical
that $\beta(1;\Q_p)=4$, and Demyanov \cite{d2}
has shown that $\beta(2;\Q_p)=8$ for all $p$.  However when $r=3$ we
only know that $\beta(3;\Q_p)=12$ for $p\ge 11$ (Schuur \cite{Sch}).
This leads us to ask the following question.
\begin{question}
Is $\beta(r;\Q_p)=4r$ for all $r$ and $p$?
\end{question}
A search for counterexamples might be worthwhile.

While Theorem \ref{odd} probably falls far short of the truth, the
hybrid method does result in a sharp bound for cubic forms.  In this
case the approach reduces to that used by Demyanov \cite{dem} in
proving that $v_3(p)=9$ for $p\not=3$.  In order to motivate our
treatment of quartic forms we reproduce our version of Demyanov's
method in section \ref{secC}, proving the following result.
\begin{theorem}\label{cub}
When $p\not=3$ we have $v_3(p)=9$.
\end{theorem}
It should be pointed out that the quasi-diagonalization aspect of our
proof of Theorem \ref{cub} only requires the solution of simultaneous
linear equations, for which we have a complete theory.  On the other
hand our treatment of Theorem \ref{odd} involves the solution of
simultaneous quadratic equations, for which our information is rather
poor.  Sharper results on the $p$-adic zeros of systems of quadratic
forms would lead to corresponding improvements in Theorem \ref{odd}.

Our approach to Theorem \ref{even} is rather different.  The method
outlined above seems hopeless for $p=2$, since we cannot exclude the
possibility that the reduction of $F(\x)$ modulo 2 is diagonal, in
which case there will only be singular zeros.  We are therefore forced
to work (essentially) with diagonal forms, as in Brauer's approach.
However we introduce a new idea which enables us to reduce
the number of variables necessary in the diagonal forms we have to produce.

A natural question is whether our results can be extended to ${\mathfrak
p}$-adic fields in general.  Our methods are in principle applicable
to these fields.  However our results rely on significant case-by-case
computer checking for forms over the residue class fields $\F_q$ with
$q<50$.  These calculations have only been carried out for prime
values of $q$.  Thus our theorems are proven only for ${\mathfrak
p}$-adic fields whose residue class field has prime order.

We introduce two points of notation which will be used throughout this
paper. Firstly, 
if $a\in\Q_p-\{0\}$ we shall use the notation $v(a)$ for the unique
integer such that $p^{-v(a)}a$ is a $p$-adic unit.  Secondly, we shall
use $\theta$ to denote the reduction map from $\Z_p$ to $\F_p$.

Finally, thanks must be recorded to the referee, who made a number of
helpful comments, and spotted a number of misprints in the original
version of this paper, as well as one significant error.

\section{Quintic Forms}\label{QF}

Our proof of Theorem \ref{Th5} is heavily based on the work of Leep
and Yeomans \cite{LY}, and our improvement stems merely from
appropriate numerical computations.  As Leep and Yeomans explain in
their introduction, they assume that $F(\x)$ is a $p$-adic quintic 
form in $n\ge 26$
variables, with only the trivial $p$-adic zero.  They then show that 
there exist $\be_1,\be_2,\be_3\in\Q_p^n$ such that if 
\[G(t_1,t_2,t_3)=\theta(F(t_1\be_1+t_2\be_2+t_3\be_3))\]
then
$G$ defines a curve with (at least) three singular points over
$\F_p$. Moreover $G$ can be taken to be absolutely irreducible if
$p\ge 7$.

Now, providing that we can find a non-singular point on $G=0$, over
$\F_p$, then this can be lifted via Hensel's Lemma to provide a
non-trivial $p$-adic solution to $F(\x)=0$.  When $p\ge 47$ Leep and
Yeomans use the Weil bound for the number of points on the curve $G=0$
to show that there is at least one non-singular point.

For each prime $p<47$ there are only finitely many forms $G$ to
consider, and one can look for a non-singular point on each of the
corresponding curves.  There cannot be three collinear singular points,
since $G$ is absolutely irreducible.  Hence we may take each of
$(1,0,0), (0,1,0)$ and $(0,0,1)$ to be singular. If $G=0$ has no
non-singular point we may then assume, after a
suitable permutation of variables, that $G$ takes one of the forms
\[G(x,y,z)=Ax^3y^2+By^3z^2+Cz^3x^2+xyzQ(x,y,z)\]
or
\[G(x,y,z)=Ax^3y^2+By^3z^2+Cz^2x^3+xyzQ(x,y,z),\]
where $Q(x,y,z)$ is quadratic.  This shows that there are essentially
9 coefficients to consider.  Allowing for the possibility of
re-scaling both the form itself and the variables, there are, in
effect just 6 degrees of freedom.

A computer calculation with forms of the above shape verifies that
whenever $17\le p< 47$ such forms always have at least one
non-singular zero, and this suffices for the theorem.  When $p=13$ the
example
\[x^3y^2+3y^3z^2+6x^3z^2+xyz(11x^2+xy+xz+6y^2+yz+4z^2)\]
shows that there need be no non-singular zero.  It seems possible that
one could tackle such cases by looking at forms $G$ in 4 variables.  However
the number of such forms appears to be too great for an exhaustive
search to succeed.

\section{Theorems \ref{odd} and \ref{even} --- Preliminaries}
\label{wooleysec}

In this section we shall explain the principles behind Wooley's
approach \cite{W} to the quasi-diagonalization procedure, and
illustrate them by verifying (\ref{med}).  We begin by introducing
some notation.  Let $S$ be any collection of $p$-adic forms in $n$
variables, comprising $r_i$ forms of degree $i$, for $1\le i\le d$.
Write $\cl{S}(n)$ for the set of such systems $S$ for which the only 
common $p$-adic zero is the trivial one.  We then define
\[V_d(r_d,r_{d-1},\ldots,r_1;p):=\max\{n:\cl{S}(n)\not=\emptyset\}.\]
Thus a system $S$ with $n>V_d(r_d,r_{d-1},\ldots,r_1;p)$ will always
have a non-trivial common zero.  We record at once the fact that
\begin{equation}\label{lin}
V_d(r_d,r_{d-1},\ldots,r_1;p)=V_d(r_d,r_{d-1},\ldots,r_2,0;p)+r_1.
\end{equation}
In addition to the above notation we shall write $\phi_d(p)$ for the largest
integer $n$ such that there is a {\em diagonal} form $F(x_1,\ldots,x_n)
\in\Q_p[x_1,\ldots,x_n]$ of degree $d$, with only the trivial $p$-adic
zero. 

Wooley's basic result \cite[Lemma 2.1]{W} is that
\begin{equation}\label{key}
V_d(r_d,r_{d-1},\ldots,r_1;p)\le\phi_d(p)+V_d(r'_d,r'_{d-1},\ldots,r_1';p),
\end{equation}
where $r'_d=r_d-1$ and
\[r'_j=\sum_{i=j}^dr_i\left(\begin{array}{c}\phi_d(p)+i-j-1\\ i-j\end{array}
\right),\;\;(1\le j\le d-1).\]
We shall present Wooley's proof in due course, since we shall need to
adapt it later.  However we begin by using (\ref{key}) to prove (\ref{med}).

From (\ref{key}) we have
\[v_4(p)=V_4(1,0,0,0;p)\le\phi+V_3(\phi,\frac{\phi(\phi+1)}{2},
\frac{\phi(\phi+1)(\phi+2)}{6};p),\]
where we have set $\phi=\phi_4(p)$ for brevity.  Moreover, writing
$\psi=\phi_3(p)$, we have
\[V_3(a,b,c;p)\le\psi+V_3(a-1,a\psi+b,a\frac{\psi(\psi+1)}{2}+b\psi+c;p),\]
whence an easy induction argument yields
\begin{eqnarray*}
\lefteqn{V_3(a,b,c;p)}\\
&\le & a\psi+V_2(\frac{a(a+1)}{2}\psi+b,
\frac{a(a+1)}{2}\frac{\psi(\psi+1)}{2}+\frac{a(a^2-1)}{3}\psi^2+ab\psi+c;p),
\end{eqnarray*}
which becomes
\begin{eqnarray}\label{med3}
V_3(a,b,c;p)&\le & a\psi+V_2(\frac{a(a+1)}{2}\psi+b,0;p)\nonumber\\
& & \hspace{3mm}\mbox{}+
\frac{a(a+1)}{2}\frac{\psi(\psi+1)}{2}+\frac{a(a^2-1)}{3}\psi^2+ab\psi+c,
\end{eqnarray}
in view of (\ref{lin}).
Finally we conclude that
\begin{eqnarray}\label{med1}
v_4(p)&\le &V_2(\frac{\phi(\phi+1)(\psi+1)}{2},0;p)\nonumber\\
&&\hspace{3mm}\mbox{}+
\phi\frac{\phi^2+3\phi+8}{6}+\psi\phi\frac{2\phi^2+3\phi+5}{4}+
\psi^2\phi\frac{4\phi^2+3\phi-1}{12}.
\end{eqnarray}

At this point we require some information about $V_2(r,0;p)$.  This
could be obtained by further applications of (\ref{key}), but in fact
rather better estimates are already available from the
literature. With this in mind we introduce the notation $\beta(r,m;K)$
for any field $K$, to denote the largest $n$ for which there are $r$
quadratic forms over $K$, in $n$ variables, having no linear space of
common zeros, defined over $K$ and having projective dimension $m$.
We also set $\beta(r;K)=\beta(r,0;K)$ which is the largest $n$ for 
which there are $r$ quadratic forms over $K$ having no common zero.
Thus $\beta(r;\Q_p)=V_2(r,0;p)$.
\begin{lemma}\label{q}
For every prime $p$ we have
\begin{enumerate}
\item[(i)] $\beta(1;\Q_p)=4$;
\item[(ii)] $\beta(2;\Q_p)=8$;
\item[(iii)] $\beta(3;\Q_p)\le 16$;
\item[(iv)] $\beta(4;\Q_p)\le 24$;
\item[(v)] $\beta(5;\Q_p)\le 40$;
\item[(vi)] $\beta(6;\Q_p)\le 56$;
\item[(vii)] $\beta(r;\Q_p)\le 2r^2-14$ for odd $r\ge 7$;
\item[(viii)] $\beta(r;\Q_p)\le 2r^2-16$ for even $r\ge 8$.
\end{enumerate}
\end{lemma}
(We are grateful to J. Zahid for
pointing out an oversight in the statement of Lemma \ref{q} in an
earlier version of this paper.) The result is a refinement of 
Corollary 2 of Dietmann \cite{Diet}, in which
we have substituted the recent result
\begin{equation}\label{uinv}
\beta(1;\Q_p(X))=8
\end{equation}
for the upper bound
\[\beta(1;\Q_p(X))\le 10\;\;\; (p\not=2)\]
of Parimala and Suresh \cite{PS}.  Wooley has proved (\ref{uinv})
using the circle method, in work to appear, while Leep \cite{Leepu},
still more recently, has given a more general result including
(\ref{uinv}) as a special case.  It should be stressed that both these
authors handle $p=2$ as well as the case of odd primes.

By using (\ref{uinv}) one may replace
\cite[(9)]{Diet} by
\begin{equation}\label{df}
\beta(r;\Q_p)\le 8+2\beta(r-2;\Q_p),
\end{equation}
which suffices for the proof of (iii)--(vi) above.  For the remaining
parts of the lemma we will use the inequalities
\begin{equation}\label{diet1}
\beta(r,m;\Q_p)\le\beta(r;\Q_p)+(r+1)m
\end{equation}
and
\[\beta(r;\Q_p)\le\beta(r-k,\beta(k;\Q_p);\Q_p)\]
of Leep \cite[Corollary 2.4, (ii)]{leep} and Martin \cite[Lemma
2]{mart}.  These yield
\[\beta(7;\Q_p)\le\beta(6,4;\Q_p)\le \beta(6;\Q_p)+28\le 84.\]
The remaining bounds (vii) and (viii)
now follow by induction from the cases $r=7$ and $r=6$ respectively, using
the bound
\begin{equation}\label{df1}
\beta(r,\Q_p)\le\beta(r-2,8;\Q_p)\le\beta(r-2;\Q_p)+8(r-1),
\end{equation}
just as in Dietmann's work.

We can do better still for $p\ge 11$, since in this case
the work of Schuur \cite{Sch} gives $\beta(3;\Q_p)=12$.  The following
result is essentially Corollary 3 of Dietmann \cite{Diet}, modified to
take account of (\ref{uinv}).
\begin{lemma}\label{qbig}
For every prime $p\ge 11$ we have 
\begin{enumerate}
\item[(i)] $\beta(3;\Q_p)=12$;
\item[(ii)] $\beta(4;\Q_p)\le 24$;
\item[(iii)] $\beta(5;\Q_p)\le 32$;
\item[(iv)] $\beta(6;\Q_p)\le 56$;
\item[(v)] $\beta(r;\Q_p)\le 2r^2-2r-12$ when $r\equiv 1\pmod{3}$ and $r\ge 7$;
\item[(vi)] $\beta(r;\Q_p)\le 2r^2-2r-8$ when $r\equiv 2\pmod{3}$ and $r\ge 8$;
\item[(vii)] $\beta(r;\Q_p)\le 2r^2-2r-8$ when 
$r\equiv 0\pmod{3}$ and $r\ge 9$.
\end{enumerate}
\end{lemma}
Here the bounds (ii) and (iv) are just parts (iv) and (vi)of Lemma 
\ref{q}, after which
part (iii) follows from (\ref{df}), as does the case $r=7$ of (v).
To obtain the case $r=9$ of part (vii) we now use
(\ref{df1}). Finally, we use  the inequalities
\[\beta(r,\Q_p)\le\beta(r-3,12;\Q_p)\le\beta(r-3;\Q_p)+12(r-2)\]
to complete the proofs of parts (v), (vi) and (vii) by induction,
starting at $r=7$, $r=5$ and $r=9$ respectively.

In order to use (\ref{med1}) we also need information about
$\phi=\phi_4(p)$ and $\psi=\phi_3(p)$.  The techniques for studying
$\phi_d(p)$ are well-known, see Davenport and Lewis \cite{DL}, for
example, so we shall merely state the following without proof.
\begin{lemma}\label{phipsi}
For $d=3$ and 4 we have
\begin{enumerate}
\item[(i)] $\phi_3(p)=3$ for $p\equiv 2\pmod{3}$;
\item[(ii)] $\phi_3(p)=6$ for $p\equiv 1\pmod{3}$; 
\item[(iii)] $\phi_3(3)=4$.
\item[(iv)] $\phi_4(p)=8$ for $p\not=2,5,13$ or $29$;
\item[(v)] $\phi_4(p)=12$ for $p=13$ or $29$;
\item[(vi)] $\phi_4(2)=15$;
\item[(vii)] $\phi_4(5)=16$.
\end{enumerate}
\end{lemma}

It is thus apparent that the worst case for (\ref{med1}) must be one
of $p=2$, $p=5$ or $p=13$.  For these values we compute that
\[v_4(2)\le \beta(480;\Q_2)+16940,\]
\[v_4(5)\le\beta(544;\Q_5)+20464\]
and
\[v_4(13)\le\beta(546;\Q_{13})+28294.\]
Moreover Lemmas \ref{q} and \ref{qbig} yield 
\[\beta(480;\Q_2)\le 460784,\;\;\;
\beta(544;\Q_5)\le 591856\;\;\;\mbox{and}\;\;\;\beta(546;\Q_{13})\le 595132,
\]
whence
\[v_4(2)\le 477724,\]
\[v_4(5)\le 612320\]
and
\begin{equation}\label{b13}
v_4(13)\le 623426.
\end{equation}
The bound (\ref{med}), stated in the introduction, now follows.

The remainder of this section will be devoted to proving (\ref{key}),
following Wooley \cite[\S 2]{W}.  We write $\phi=\phi_d(p)$ for short,
and suppose that 
\[n>\phi+V_d(r'_d,r'_{d-1},\ldots,r_1';p). \]
Let 
\[\tilde{r}_i=\left\{\begin{array}{cc} r_d-1,&\;\; i=d,\\ r_i,&\;\;
    i<d, \end{array}\right.\]
and suppose our system $S$ consists of a form $F$ of degree $d$ along
with forms $G_{ij}$ of degree $i$ for $1\le j\le \tilde{r}_i$ and
$1\le i\le d$.  By using induction on $k$ we shall find linearly
independent vectors $\be_1,\ldots,\be_k\in\Q_p^n$ such that
$F(t_1\be_1+\ldots+t_k\be_k)$ is a diagonal form in $t_1,\ldots,t_k$,
and for which each form $G_{ij}(t_1\be_1+\ldots+t_k\be_k)$ vanishes
identically. If we can do this for $k=1+\phi$ then an appropriate
choice of the $t_i$ will make every form in the system vanish, as
required.

Since $n>V_d(r'_d,r'_{d-1},\ldots,r_1';p)\ge V_d(\tilde{r}_d,\ldots,
\tilde{r}_1;p)$ we can find a non-zero vector $\be_1$ at which every
form $G_{ij}$ vanishes.  This is enough to establish the base case
$k=1$ for the induction.  Now suppose that $k\le\phi$, and that we
have found a suitable set of vectors $\be_1,\ldots,\be_k$.  Let
$T\subseteq\Q_p^n$ be the space spanned by $\be_1,\ldots,\be_k$, and
take $U$ to be any direct complement of $T$, so that $T\oplus
U=\Q_n^p$. We shall insist that $\be_{k+1}\in U-\{\b{0}\}$, so that 
$\be_1,\ldots,\be_k,\be_{k+1}$ will automatically be linearly
independent.  We also note that
\begin{equation}\label{dimU}
{\rm dim}(U)=n-k\ge n-\phi>V_d(r'_d,r'_{d-1},\ldots,r_1';p).
\end{equation}
For each multi-degree vector $\b{u}=(u_1,\ldots,u_k)$, where the $u_i$
are non-negative integers, we write 
$|\b{u}|=u_1+\ldots+u_k$.
We then proceed to define forms $F^{(\lb{u})}$ by the
expansion
\[F(t_1\be_1+\ldots+t_k\be_k+t\x)=
\sum_{|\lb{u}|\le d}\b{t}^{\lb{u}}t^{d-|\lb{u}|}F^{(\lb{u})}(\x),\]
where we have written
\[\b{t}^{\lb{u}}=\prod_{i=1}^k t_i^{u_i}.\]
Thus $F^{(\lb{u})}(\x)$ will be a form of degree $d-|\b{u}|$.
Similarly we define forms $G_{ij}^{(\lb{u})}$ by writing
\[G_{ij}(t_1\be_1+\ldots+t_k\be_k+t\x)=
\sum_{|\lb{u}|\le i}\b{t}^{\lb{u}}t^{i-|\lb{u}|}G_{ij}^{(\lb{u})}(\x),\]
so that $G_{ij}^{(\lb{u})}$ has degree $i-|\b{u}|$.  We now see that
$\be_{k+1}=\x\in U-\{\b{0}\}$ will be an admissible choice providing
that
\begin{equation}\label{f}
F^{(\lb{u})}(\x)=0\;\;\;\mbox{for all}\;\b{u}\not=\b{0}
\end{equation}
and
\begin{equation}\label{gij}
G_{ij}^{(\lb{u})}(\x)=0\;\;\;\mbox{for all}\;\b{u}.
\end{equation}
Thus $\x$ must be a common zero of a new system of forms $S'$, say.
It remains to check how many forms there are of each degree.  The only
forms of degree $d$ arise from (\ref{gij}) with $i=d$ and
$\b{u}=\b{0}$.  There are therefore $r_d-1=r_d'$ such forms. In
general the number of vectors $\b{u}$ with $|\b{u}|=u$ is
\[\left(\begin{array}{c}u+k-1\\ u\end{array}\right).\]
Thus, for $m<d$, we get
\[\left(\begin{array}{c}d-m+k-1\\ d-m\end{array}\right)\le
\left(\begin{array}{c}\phi+d-m-1\\ d-m\end{array}\right)\]
forms of degree $m$ from (\ref{f}), and
\[\left(\begin{array}{c}i-m+k-1\\ i-m\end{array}\right)\le
\left(\begin{array}{c}\phi+i-m-1\\ i-m\end{array}\right)\]
such forms from (\ref{gij}), for each $i$ and $j$.  The system $S'$
therefore consists of at most $r_m'$ forms of degree $m$, for $1\le m\le d$.
In view of (\ref{dimU}) there is therefore a suitable common solution
$\x$, which completes our induction step.

The above is the argument as Wooley presents it, however we observe
that a small saving can be made by requiring only that
$\be_{k+1}\in\Q_p^n-\{\b{0}\}$, rather than $\be_{k+1}\in
U-\{\b{0}\}$. With this change it is no longer immediate that
$\be_1,\ldots,\be_k,\be_{k+1}$ are linearly independent.  However if
there is a dependence relation we may write it as $\be_{k+1}=
\sum_{i=1}^k c_i\be_i$, since our induction assumption shows that
$\be_1,\ldots,\be_k$ are linearly independent. We now choose
$\x=\be_{k+1}$ to be a non-zero vector satisfying (\ref{f}) and
(\ref{gij}) as before, whence we will have
\begin{equation}\label{n1}
F(t_1\be_1+\ldots+t_k\be_k+t_{k+1}\be_{k+1})=\sum_{i=1}^{k+1}A_it_i^{d}
\end{equation}
and
\begin{equation}\label{n2}
G_{ij}(t_1\be_1+\ldots+t_k\be+t_{k+1}\be_{k+1})=0
\end{equation}
identically in $t_1,\ldots,t_{k+1}$.  On substituting for $\be_{k+1}$
in the first of these relations we would find that
\begin{eqnarray*}
\sum_{i=1}^{k+1}A_it_i^{d}&=&F(t_1\be_1+\ldots+t_k\be_k+t_{k+1}\be_{k+1})\\
&=&F(\{t_1+c_1t_{k+1}\}\be_1+\ldots+\{t_k+c_k t_{k+1}\}\be_k)\\
&=&\sum_{i=1}^{k}A_i(t_i+c_it_{k+1})^{d},
\end{eqnarray*}
identically in $t_1,\ldots,t_{k+1}$.  Thus we must have $A_ic_i=0$ for
each $i\le k$.  Since $\be_{k+1}\not=\b{0}$ there must be at least one
non-zero value of $c_i$, so that $A_i=0$ for some index $i=i_0$,
say.  However, it then follows from (\ref{n1}) and (\ref{n2}) that
$\be_{i_0}$ is a common zero of the system $S$.  Thus, either
$\be_1,\ldots,\be_{k+1}$ are indeed linearly independent, or we have a
suitable common zero for our system.

It follows that Wooley's estimate (\ref{key}) can be replaced by
\[V_d(r_d,r_{d-1},\ldots,r_1;p)\le V_d(r'_d,r'_{d-1},\ldots,r_1';p).\]
As a result we may replace (\ref{med3}) by
\begin{eqnarray}\label{med4}
V_3(a,b,c;p)&\le & V_2(\frac{a(a+1)}{2}\psi+b,0;p)\nonumber\\
& & \hspace{3mm}\mbox{}+
\frac{a(a+1)}{2}\frac{\psi(\psi+1)}{2}+\frac{a(a^2-1)}{3}\psi^2+ab\psi+c,
\end{eqnarray}
and (\ref{med1}) by
\begin{eqnarray*}
v_4(p)&\le &V_2(\frac{\phi(\phi+1)(\psi+1)}{2},0;p)\\
&&\hspace{3mm}\mbox{}+
\phi\frac{\phi^2+3\phi+2}{6}+\psi\phi\frac{2\phi^2+3\phi+1}{4}+
\psi^2\phi\frac{4\phi^2+3\phi-1}{12}.
\end{eqnarray*}

A further small saving can be obtained by observing that
\[V_3(1,b,0;p)\le \beta(b,9;\Q_p).\]
To prove this, suppose we are given a system $S$ consisting of a 
cubic form $C$ and
quadratic forms $Q_1,\ldots,Q_b$.  Suppose further that we have 
sufficient variables that the quadratic forms have a linear space $L$ 
of common zeros, where $L$ has projective dimension 9.  Then $C$ will
vanish on $L$ since we may take $v_3=9$, and hence the system $S$ has a
common zero.  Now (\ref{diet1}) yields
\[V_3(1,b,0;p)\le \beta(b;\Q_p)+9(b+1).\]
If we use this to start the induction, we replace (\ref{med4}) by
\[V_3(a,b,c;p)\le  V_3(1,b',c';p)\]
with
\[b'=\frac{(a-1)(a+2)}{2}\psi+b\]
and
\[c'=\frac{(a-1)(a+2)}{2}\frac{\psi(\psi+1)}{2}+
\frac{(a-1)(a-2)(2a+3)}{6}\psi^2+(a-1)b\psi+c.\]
Hence
\begin{equation}\label{med5}
V_3(a,b,c;p)\le\beta(\frac{(a-1)(a+2)}{2}\psi+b;\Q_p)+c'',
\end{equation}
with
\begin{eqnarray*}
c''&=&9(\frac{(a-1)(a+2)}{2}\psi+b+1)+
\frac{(a-1)(a+2)}{2}\frac{\psi(\psi+1)}{2}\\
&&\hspace{3mm}\mbox{}+
\frac{(a-1)(a-2)(2a+3)}{6}\psi^2+(a-1)b\psi+c.
\end{eqnarray*}
These minor variants result in a rather small overall improvement.
Thus we may replace (\ref{b13}) by
\[v_4(13)\le 611930\]
for example.

\section{Cubic Forms}\label{secC}

In this section we shall develop our hybrid approach to Artin's
problem, and illustrate it in its simplest setting by proving Theorem
\ref{cub}. We shall argue by contradiction, and so we
suppose that $F(\x)\in\Q_p[\x]$ is a form of degree 3, in 10 variables, 
with only
the trivial $p$-adic zero.  Our overall strategy will be to seek
linearly independent vectors $\be_1,\be_2,\be_3\in\Q_p^{10}$ such that,
for an appropriate $r\in\Z$, the form
$p^{-r}F(x\be_1+y\be_2+z\be_3)$ has coefficients in $\Z_p$, and  such that
$\theta(p^{-r}F(x\be_1+y\be_2+z\be_3))$ has at least one non-singular zero.
In particular it will follow by
Hensel's Lemma that $p^{-r}F(x\be_1+y\be_2+z\be_3)$ has a non-trivial
$p$-adic zero, and hence that $F(\x)$ similarly has a non-trivial zero.

When
$\x\in\Q_p^{10}-\{\b{0}\}$ we shall say that $\x$ has ``level $r$'',
where $0\le r\le 2$, if $v(F(\x))\equiv r\pmod{3}$.  Since we are
assuming that $F(\x)\not=0$ for such $\x$, this concept is
well-defined.  For any set
\[S=\{\be_1,\ldots,\be_m\}\subset\Q_p^{10}-\{\b{0}\}\]
we say that $S$ is
``admissible'' if
\begin{enumerate}
\item[(i)] $0\le v(F(\be_i))\le 2$ for $1\le i\le m$.
\item[(ii)] For each level $r$ there are at most two vectors
  $\be_i$ of level $r$.
\item[(iii)] If $\be_i$ and $\be_j$ have the same level, with $i<j$,
  then 
\[F(x\be_i+y\be_j)=Ax^3+Bxy^2+Cy^3\]
for certain $A,B,C\in\Q_p$  depending on $i$ and $j$.
\end{enumerate}
It is clear that if $\be\not=\b{0}$ then the singleton set $S=\{p^{-k}\be\}$
is admissible for some $k$.  
Moreover any admissible set has cardinality at most 6,
by (ii).  If $\be_i$ and $\be_j$ have the same level, they must be
linearly independent, by the following result.
\begin{lemma}\label{li}
Let $F(\x)\in\Q_p[x_1,\ldots,x_n]$ be a form of degree $d$, having
only the trivial zero in $\Q_p^n$. Let
$\be_1,\ldots,\be_k$ be linearly independent vectors in $\Q_p^n$, and
suppose we have a non-zero vector $\be\in\Q_p^n$ such that the form
\[F_0(t_1,\ldots,t_k,t):=F(t_1\be_1+\ldots+t_k\be_k+t\be)\]
in the indeterminates $t_1,\ldots,t_k$ and $t$,
contains no terms of degree one in $t$.  Then the set
$\{\be_1,\dots,\be_k,\be\}$ is linearly independent.
\end{lemma}
In order not to interrupt our discussion of cubic forms we postpone the
proof of this until the end of the present section.

Before proceeding further we note that
if $\be_i$ and $\be_j$ both have level $r$, say, then
$p^{-r}F(x\be_i+y\be_j)$ must have coefficients in $\Z_p$.  This
follows from our next result.
\begin{lemma}\label{mcoe}
Let $f(x,y)=ax^d+bxy^{d-1}+cy^d\in\Q_p[x,y]$ and suppose that
$a,c\in\Z_p$, but that $b\not\in\Z_p$.  Then there exist
$\alpha,\beta\in\Q_p$, not both zero, for which $f(\alpha,\beta)=0$.
\end{lemma}
This too we will prove at the end of the section.

We now assume that we have an admissible set $S$ of maximal
size. We seek one further non-zero vector $\be\in\Q_p^{10}$, satisfying
certain further constraints, which will correspond to the
quasi-diagonalization step.  There are constraints for each of the
three levels $r=0,1,2$, which we now describe.  
If the set $S$ has no elements of level $r$
there will be no corresponding constraints.  If $S$ has exactly one
element, $\be_i$ say, of level $r$ we write
\[F(x\be_i+y\be)=x^3F(\be_i)+x^2yL_i(\be)+xy^2Q_i(\be)+y^3F(\be),\]
where $L_i(\be)$ is a linear form in $\be$, depending on $\be_i$, and
$Q_i(\be)$ is similarly a quadratic form in $\be$, depending on
$\be_i$.  In this case we shall impose on $\be$ the single linear
constraint $L_i(\be)=0$.

When $S$ has two elements $\be_i,\be_j$ of level $r$ we write
\begin{eqnarray*}
F(x\be_i+y\be_j+z\be)&=&F(x\be_i+y\be_j)+
\{x^2L_i(\be)+xyL_{ij}(\be)+y^2L_j(\be)\}z\\
&&\hspace{1cm}\mbox{}+\{xQ_i(\be)+yQ_j(\be)\}z^2+
F(\be)z^3,
\end{eqnarray*}
where $L_i,L_{ij},L_j$ are linear forms and $Q_i,Q_j$ are quadratic
forms.  In this case we impose the three linear constraints
$L_i(\be)=L_{ij}(\be)=L_j(\be)=0$. 

Thus $\be$ has to satisfy at most 9 linear constraints, so that we may
indeed find a suitable $\be\in\Q_p^{10}-\{\b{0}\}$.  We now recall
that $S$ was chosen to be maximal.  By construction we therefore
see that if $\be$ is of level $r$ then there must have been two
vectors $\be_i,\be_j$ in $S$ which also have level $r$.  We take
$i<j$, and multiply $\be$ by an appropriate power of $p$ so that
$v(F(\be))=r$. After changing notation slightly
from (iii) above we may then write 
\[p^{-r}F(x\be_i+y\be_j+z\be)=Ax^3+Bxy^2+Cy^3+(Dx+Ey)z^2+Fz^3,\]
where $A,C,F$ are $p$-adic units.  We noted earlier that $B$ must be a
$p$-adic integer. Similarly, taking $y=0$, Lemma \ref{mcoe} shows that
$D$ is a $p$-adic integer.  Setting $x=0$ yields the same
conclusion for $E$. Moreover $\be_i,\be_j$ and $\be$ must be linearly
independent by Lemma \ref{li}.
We now call on the following lemma.
\begin{lemma}\label{cl}
Let $p\not=3$ and suppose that
\[f(x,y,z)=ax^3+bxy^2+cy^3+(dx+ey)z^2+fz^3\in\F_p[x,y,z],\]
with $acf\not=0$.  Then $f$ has at least one non-singular zero over
$\F_p$.
\end{lemma}
If we use this in conjunction with Hensel's Lemma we find that
$F(x\be_i+y\be_j+z\be)$ has a non-trivial $p$-adic zero.  Thus $F(\x)$
also has a non-trivial zero, which completes the proof of Theorem
\ref{cub}.

It remains to prove Lemmas \ref{li}, \ref{mcoe} and \ref{cl}, and 
we begin with the first of these. We suppose for a contradiction that
$\be=a_1\be_1+\ldots+a_k\be_k$.  We would then have
\begin{eqnarray*}
(1+t)^dF(\be)&=&F((1+t)\be)\\
&=&F(a_1\be_1+\ldots+a_k\be_k+t\be).
\end{eqnarray*}
By our hypothesis, the final expression contains no linear term in
$t$, while the first expression contains the term $dF(\be)t$.  Thus we
must have $F(\be)=0$, contradicting the assumption that $F(\x)$ has
only the trivial zero.

Next we examine Lemma \ref{mcoe}.  Suppose that $v(b)=s<0$.  Then
$p^{-s}f(x,y)\equiv b'xy^{d-1}\pmod{p}$, where $b'$ is a $p$-adic
unit. Thus $\theta(p^{-s}f(x,y))$ has a
non-singular zero at $(0,1)$, from which Hensel's Lemma produces the
required solution $f(\alpha,\beta)=0$ in $\Q_p$.

Finally we prove Lemma \ref{cl}. Suppose firstly that $f$ is 
absolutely irreducible.  Write $N$ for the number of points over
$\F_p$, lying on the projective curve $f=0$. By the Weil bound in the
form given by Leep and Yeomans
\cite[Corollary 1]{LYc}, we have
\[|N-(p+1)|\le 2g\sqrt{p}+1-g,\]
where $g=0$ or 1.  Since there is no singular point when the genus $g$
is 1, and one singular point when $g=0$, we conclude that there is
always at least one non-singular point, as required.

If $f$ factors over $\overline{\F_p}$ it must have a linear factor,
$z-L(x,y)$, say.  Then $f(x,y,L(x,y))$ will vanish identically, so
that $L(x,y)^2$ divides $ax^3+bxy^2+cy^3$.  Now $ax^3+bxy^2+cy^3$
cannot be a multiple of $L(x,y)^3$, since it has no term in $x^2y$
and $p\not=3$.
Hence it must have a linear factor, $L'(x,y)$ say, of multiplicity one.
Moreover if we multiply $L$ and $L'$ by appropriate constants it is clear that
they must be defined over $\F_p$.  Thus $ax^3+bx+c$ has a root, $u$
say, of multiplicity one and lying in $\F_p$.  It then follows that
$(u,1,0)$ is a non-singular zero of $f$.

\section{Theorem \ref{odd} --- A Preliminary Lemma}\label{largeodd}

In the next two sections we shall consider Theorem \ref{odd} for
$p\not=5$.
We begin by proving the following key result.  It will be convenient
to say that two forms $f(x_1,\ldots,x_m)$ and $g(x_1,\ldots,x_m)$ over
a field $F$ are ``similar'' if there are non-zero elements
$a,a_1,\ldots,a_m\in F$ such that
\[f(x_1,\ldots,x_m)=ag(a_1x_1,\ldots,a_mx_m).\]
\begin{lemma}\label{mykey}
Let $p\not\in\{2,5\}$ be a prime, and let
\[f(x,y)=Ax^4+Bxy^3+Cy^4\in\F_p[x,y]\]
be a binary quartic form with $AC\not=0$.
Then there exists a quadratic form $q(x,y)\in\F_p[x,y]$ with the
following properties.
\begin{enumerate}
\item[(i)] $q(x,y)$ factors over $\F_p$ into distinct linear factors.
\item[(ii)] For any $D,E,F,G\in\F_p$ with $G\not=0$, if the form 
\begin{equation}\label{gf}
g(x,y,z):=f(x,y)+Dq(x,y)z^2+Exz^3+Fyz^3+Gz^4
\end{equation}
does not have any non-singular zero over $\F_p$ then either
$p\in\{5,13\}$ and $g$ is diagonal, or $p\equiv 5$ or $7\pmod{8}$ and
$g$ is similar to 
\begin{equation}\label{bf}
x^4-4xy^3+3y^4+4H(x-y)yz^2+2H^2z^4
\end{equation}
for some $H\in\F_p-\{0\}$.
\end{enumerate}
\end{lemma}
In proving Theorem \ref{odd} we will use the form $g(x,y,z)$ in place
of a diagonal ternary quartic form.  Producing such forms $g$ from the
original quartic $F(\x)$ will require distinctly fewer variables than
would be needed to produce a diagonal form.  For the proof of Lemma
\ref{mykey} we consider four cases.

{\bf Case 1.}
This is the case in which $p\le 31$, so that $p=3,7,11,17,19,23$ or 31.  
For these primes
the theorem is proved by a computer search over all forms
$f$, in every case finding an acceptable quadratic $q$.  Thus for the
remainder of our treatment we shall assume that $p\ge 37$.

{\bf Case 2.}
Suppose next that $f(x,1)$ has a root $\xi\in\F_p$, of multiplicity one.
Then $(\xi,1)$ will be a non-singular zero of $f$, so that
$(\xi,1,0)$ will be a non-singular zero of $g$ irrespective of
the choice of $q$ or of $D,E,F$ and $G$.  Hence in this case we may
choose $q(x,y)=x(x+y)$, for example.

{\bf Case 3.}
The main case is that in which $f(x,1)$ does not have a root in
$\F_p$, and does not have a repeated root in $\overline{\F_p}$.
We begin by observing
that there must be at least one value $\alpha\in\F_p$ for which
\begin{equation}\label{ne}
\left(\frac{f(\alpha,1)}{p}\right)\not=\left(\frac{C}{p}\right),
\end{equation}
for if not, the equation $f(X,1)=CY^2$ would have exactly $2p$
solutions over $\F_p$.  Since $f$ has no repeated factor this would
contradict the Weil bound, since $|2p-(p+1)|>2p^{1/2}$ for $p\ge 37$.
We fix an $\alpha$ for which (\ref{ne}) holds, and note that
$\alpha\not=0$, since $f(0,1)=C$.  We then define
$q(x,y)=x(x-\alpha y)$, which clearly satisfies part (i) of the
lemma. It therefore remains to verify
part (ii).  

We begin by showing that the form $g(x,y,z)$ must be
absolutely irreducible. Our first step is
to demonstrate that $g$ cannot have quadratic factors over
$\overline{\F_p}$.  Suppose
\begin{eqnarray*}
g(x,y,z)&=&Gz^4+z^3(Ex+Fy)+Dz^2q(x,y)+f(x,y)\\
&=& G(z^2+zL_1(x,y)+Q_1(x,y))(z^2+zL_2(x,y)+Q_2(x,y)),
\end{eqnarray*}
with $L_1,L_2$ linear and $Q_1,Q_2$ quadratic.  Then
\begin{equation}\label{c1}
L_1L_2+Q_1+Q_2=DG^{-1}q,
\end{equation}
\begin{equation}\label{c2}
L_1Q_2+L_2Q_1=0,
\end{equation}
and
\begin{equation}\label{c3}
Q_1Q_2=G^{-1}f.
\end{equation}
Now, since $f$ does not have a repeated factor over $\overline{\F_p}$,
it follows from (\ref{c3}) that $Q_1$ and $Q_2$ are coprime.  We may
then deduce from (\ref{c2}) that $L_1=L_2=0$.  Hence in order to
solve (\ref{c1}) we set $Q_1=DG^{-1}q/2+R$ and $Q_2=DG^{-1}q/2-R$, with $R\in
\overline{\F_p}[x,y]$. Thus (\ref{c3}) produces $G^{-1}f=
D^2G^{-2}q^2/4-R^2$, so
that in fact $R$ takes the shape $R=k^{1/2}S$ with $k\in\F_p$ and
$S\in\F_p[x,y]$. Now, if we set $(x,y)=(0,1)$ in the relation 
\[G^{-1}f(x,y)=D^2G^{-2}q(x,y)^2/4-kS(x,y)^2,\]
and recall that $q(x,y)=x(x-\alpha y)$, we find that $C=-kGS(0,1)^2$, whence
\[\left(\frac{C}{p}\right)=\left(\frac{-kG}{p}\right).\]
On the other hand, if we take $(x,y)=(\alpha,1)$ we obtain
\[\left(\frac{f(\alpha,1)}{p}\right)=\left(\frac{-kG}{p}\right).\]
(Note that one cannot have $f(\alpha,1)=0$, since $f(x,1)$ has no
roots in $\F_p$ in Case 3.)  We have thus obtained a contradiction to
(\ref{ne}), showing that $g(x,y,z)$ cannot factor into two quadratics.

It now readily follows that $g(x,y,z)$ must be absolutely
irreducible.  For otherwise it must factor into a linear form and a
cubic form, both defined over $\F_p$.  This would imply that $f(x,y)$
also has a linear factor over $\F_p$, which is contrary to the
hypotheses for Case 3.  Now suppose that the projective curve defined
over $\F_p$ by $g(x,y,z)=0$ has genus $g$, and $N$ points over $\F_p$,
of which $S$ are singular.  Then, according to Leep and Yeomans
\cite[Corollary 1 \& Lemma 1]{LYc}, we have
\[|N-(p+1)|\le 2g\sqrt{p}+3-g\]
and $0\le g\le 3-S$.  If all the points on the curve were singular we
would have $N=S$.  If $S=0$ this yields $g\le 3$ and $p+1\le
6\sqrt{p}$, which is impossible for $p\ge 37$.  On the other hand if 
$1\le S=N\le 3$ we have $g\le 2$ and 
\[|S-(p+1)|\le 4\sqrt{p}+3.  \]
This
would lead to $p-2\le 4\sqrt{p}+3$, which is also impossible for $p\ge 37$.
Hence in either case we find that $N$ cannot be equal to $S$.  Thus
the curve must have at least one non-singular point, which suffices for
(ii) of the lemma.

{\bf Case 4.}
The remaining case is that in which $f(x,1)$ does not have a root of
multiplicity one in
$\F_p$, but has a repeated root, $\rho$ say, in $\overline{\F_p}$.
Since $f$ has no term in $x^3y$ or
$x^2y^2$ it must take the shape
\begin{equation}\label{fact}
f(x,y)=Ax^4+Bxy^3+Cy^4=A(x-\rho y)^2(x^2+2\rho xy+3\rho^2 y^2),
\end{equation}
whence $B=-4A\rho^3$ and $C=3A\rho^4 $.  Since $AC\not=0$ we see that
$\rho$ and $B$ are nonzero, and hence that $\rho=-4C/(3B)\in\F_p$.  We
can therefore re-scale the form $f$ and the variable $y$ so as to
assume that 
\[f(x,y)=x^4-4xy^3+3y^4=(x-y)^2(x^2+2xy+3y^2).\]
It is clear
that $f(x,y)$ cannot have $x-y$ as a factor of multiplicity 3 or
more, since $x^2+2x+3$ cannot vanish at $x=1$.  Moreover
$x^2+2x+3$ cannot be a square, and it has no roots in $\F_p$, since 
we are not in Case 2.  It follows that $-2$ is not a quadratic residue
of $p$, so that Case 4 can arise only when $p\equiv 5$ or $7\pmod{8}$.

We shall take $q(x,y)=(x-y)y$, which clearly satisfies (i) of the lemma.
We proceed to demonstrate that it also satisfies (ii).
As in Case 3 we shall show that the form $g(x,y,z)$ must be
absolutely irreducible, unless it is similar to a form of the type
described.  Again we begin by considering quadratic factors
over $\overline{\F_p}$.  Thus we examine the
conditions (\ref{c1}), (\ref{c2}) and (\ref{c3}) as before.
Since $f(x,y)$ has a factor $x-y$, we have
$x-y|Q_1(x,y)$, say, by (\ref{c3}). We also have $x-y|q(x,y)$ by
construction.  Thus (\ref{c1}) and (\ref{c2}) yield
\[x-y|L_1L_2+Q_2,\;\;\;\mbox{and}\;\;\; x-y|L_1Q_2,\]
whence $x-y|Q_2$.  Moreover $x-y$ must divide at least one
of $L_1$ and $L_2$.  Indeed since $(x-y)^3$ does not divide
$Q_1Q_2=G^{-1}f$ it follows from (\ref{c2}) that $x-y$ divides both
$L_1$ and $L_2$.  The forms $Q_1$ and $Q_2$ cannot be proportional,
since $f(x,y)=(x-y)^2(x^2+2xy+3y^2)$ is not a square over $\overline{\F_p}$.
It therefore follows that $L_1$ and $L_2$ both vanish.

We now have $Q_1Q_2=G^{-1}f$ and $Q_1+Q_2=DG^{-1}(x-y)y$.  Thus 
\[D^2G^{-2}(x-y)^2y^2-4G^{-1}f=(Q_1-Q_2)^2\]
is a square over $\overline{\F_p}$, and hence so is
$D^2G^{-2}y^2-4G^{-1}(x^2+2xy+3y^2)$.  This latter expression is
therefore of the form $a(x+by)^2$, in which we must have $a=-4G^{-1}$
and $b=1$ in order for the coefficients of $x^2$ and $xy$ to match.
Equating the coefficients of $y^2$ then yields $D^2=8G$, whence $g$
has the shape described in the lemma.

We now see that if $g(x,y,z)$ is not absolutely irreducible, and is
not of the exceptional shape described in the lemma, then it
must factor as the product of a linear form and a
cubic form, both defined over $\F_p$. If we write $L(x,y,z)$
for the linear form then we have $L(x,y,0)|f(x,y)$. By the
hypotheses of Case 4, the only root of $f(x,1)$ in $\F_p$ is $x=1$,
whence we may take $L(x,y,0)=x-y$.  We may therefore write
$L(x,y,z)=x-y-\pi z$, where $\pi\not=0$, in view of the fact 
that $G\not=0$. Since $L|g$, the form 
\begin{eqnarray*}
g(x,y,\pi^{-1}(x-y))&=&f(x,y)+Dq(x,y)\pi^{-2}(x-y)^2\\
&& \hspace{1cm}\mbox{}+
(Ex+Fy)\pi^{-3}(x-y)^3+G\pi^{-4}(x-y)^4
\end{eqnarray*}
must vanish identically.  This however is impossible because 
$x-y|q(x,y)$ while $(x-y)^3\nmid f(x,y)$.  

Thus $g(x,y,z)$ is absolutely irreducible, and we may now prove
(ii) as in Case 3.  This completes the argument for Lemma \ref{mykey}.

\section{Theorem \ref{odd} --- $p\not=5$}\label{large}

We turn now to the proof of Theorem \ref{odd} for primes $p\not=5$.
Our goal will be to prove the following
estimate.
\begin{lemma}\label{Ll}
For primes $p\not\in\{2,5\}$ we have
\[v_4(p)\le 16+\beta(8;\Q_p).\]
\end{lemma}
On combining this with the case $r=8$ of Lemma \ref{q} or \ref{qbig}
as appropriate, we obtain the corresponding result in Theorem \ref{odd}.

It is of interest to note that there is an easy lower bound for $v_4(p)$ of
a rather similar flavour.
\begin{lemma}\label{lb}
For every prime $p$ we have
\[v_4(p)\ge \beta(4;\Q_p).\]
\end{lemma}
To prove this we take a set of $p$-adic quadratic forms
$q_i(x_1,\ldots,x_m)$ for $1\le i\le 4$ having no common $p$-adic zero
apart from the trivial one, and in which $m$ has its maximal value
$m=\beta(4;\Q_p)$. Then if $Q(y_1,\ldots,y_4)$ is anisotropic over
$\Q_p$ the quartic form
\[F(\x)=Q(q_1(\x),q_2(\x),q_3(\x),q_4(\x))\]
will have no non-trivial zero, and the lemma follows.

To prove Lemma \ref{Ll} we shall follow the method given previously for
Theorem \ref{cub}, but with an additional twist, to cover the
exceptional cases in Lemma \ref{mykey}.  We argue by contradiction, and so we
suppose that $F(\x)\in\Q_p[\x]$ is a form of degree 4 with only
the trivial $p$-adic zero.  Our overall strategy will be to seek
linearly independent vectors $\be_1,\be_2,\be_3\in\Q_p^n$ such that,
for an appropriate $r\in\Z$, the forms $p^{-r}F(x\be_1+y\be_2)$ and
$p^{-r}F(x\be_1+y\be_2+z\be_3)$ have coefficients in $\Z_p$, and their
reductions modulo $p$ are of the shape $f(x,y)$ and $g(x,y,z)$
described in Lemma \ref{mykey}.  In particular, unless we are in an
exceptional case, it will follow by
Hensel's Lemma that 
\[p^{-r}F(x\be_1+y\be_2+z\be_3)\]
has a non-trivial
$p$-adic zero, and hence that $F(\x)$ similarly has a non-trivial zero.

As before, when
$\x\in\Q_p^n-\{\b{0}\}$ we shall say that $\x$ has ``level $r$''
if $v(F(\x))\equiv r\pmod{4}$ with $0\le r\le 3$.  Since we are
assuming that $F(\x)\not=0$ for such $\x$, this concept is
well-defined.  For any set
$S=\{\be_1,\ldots,\be_m\}\subset\Q_p^n-\{\b{0}\}$ we say that $S$ is
``admissible'' if
\begin{enumerate}
\item[(i)] $0\le v(F(\be_i))\le 3$ for $1\le i\le m$.
\item[(ii)] For each level $r$ there are at most two vectors
  $\be_i$ of level $r$.
\item[(iii)] The set of all vectors $\be_i$ of a 
given level is linearly independent.
\item[(iv)] If $\be_i$ and $\be_j$ are both of level $r$, with $i<j$,
  then the form $p^{-r}F(x\be_i+y\be_j)$ has coefficients in
  $\Z_p$, and $\theta(p^{-r}F(x\be_i+y\be_j))=Ax^4+Bxy^3+Cy^4$ for 
certain $A,B,C\in\F_p$
  depending on $i$ and $j$.
\end{enumerate}
This definition is not quite the obvious modification of that given in
\S \ref{secC}.   We
shall say that a level $r$ for which there are exactly two
vectors $\be_i$ and $\be_j$ is ``suitable'', unless 
$p\equiv 5$ or $7\pmod{8}$ and $Ax^4+Bxy^3+Cy^4$ is similar to 
$x^4-4xy^3+3y^4$.  Moreover, we shall say that a level for which 
there are exactly two vectors $\be_i$ and $\be_j$ is
``acceptable'' unless $p\in\{5,13\}$ and $B=0$. 

Of all admissible sets $S$, we consider those of maximal size.  Of all 
such sets we examine those with as few
unsuitable levels as possible, and from these we select one with as
few unacceptable levels as possible. As in \S \ref{secC} we proceed to
produce a further
non-zero vector $\be$ satisfying certain constraints, which we now
describe. 

If the set $S$ has no elements of level $r$
there will be no corresponding constraints.  If $S$ has exactly one
element, $\be_i$ say, of level $r$ we write
\[F(x\be_i+y\be)=x^4F(\be_i)+x^3yL_i(\be)+x^2y^2Q_i(\be)
+xy^3C_i(\be)+y^4F(\be),\]
where $L_i,Q_i,C_i$ are forms in $\be$, depending on $\be_i$, of
degrees 1, 2 and 3 respectively.  In this case we shall impose on
$\be$ the constraints $L_i(\be)=Q_i(\be)=0$.

When $S$ has two elements $\be_i,\be_j$ of level $r$ we have more work
to do.  We take $f(x,y)=\theta(p^{-r}F(x\be_i+y\be_j))$, so that 
$f(x,y)$ satisfies the hypotheses of
Lemma \ref{mykey}.  The lemma then produces a quadratic form
$q(x,y)\in\F_p[x,y]$, which will depend on $i$ and $j$.  Let
$Q(x,y)\in\Z_p[x,y]$ be any lift of $q(x,y)$. Since $q(x,y)$ does not
vanish identically, the coefficients of $Q(x,y)$ are $p$-adic
integers, at least one of which is a $p$-adic unit.  We write 
$Q(x,y)=M_{11}x^2+M_{12}xy+M_{13}y^2$.  Then there is a $3\times 3$
unimodular matrix $M=(M_{ij})$ with entries in $\Z_p$.  We define
quadratic forms $Q'(x,y),Q''(x,y)\in\Z_p[x,y]$ by the equation
\[M\left(\begin{array}{c} x^2\\ xy\\ y^2 \end{array}\right)=
\left(\begin{array}{c} Q(x,y)\\ Q'(x,y)\\ Q''(x,y)
  \end{array}\right).\]
Thus if $N=M^{-1}$ then $N$ has $p$-adic integer entries and
\begin{equation}\label{inv}
N\left(\begin{array}{c} Q(x,y)\\ Q'(x,y)\\ Q''(x,y)\end{array}\right)
=\left(\begin{array}{c} x^2\\ xy\\ y^2 \end{array}\right).
\end{equation}
We now write
\begin{eqnarray}\label{5a}
F(x\be_i+y\be_j+z\be)&=&F(x\be_i+y\be_j)+F_3(x,y;\be)z+F_2(x,y;\be)z^2
\nonumber\\
&&\hspace{2cm}\mbox{}+F_1(x,y;\be)z^3+F(\be)z^4,
\end{eqnarray}
where each $F_i(x,y;\be)$ is bi-homogeneous, of degree $i$ in $(x,y)$
and of degree $4-i$ in $\be$.  In particular we have
\begin{equation}\label{5b}
F_3(x,y;\be)=x^3L_1(\be)+x^2yL_2(\be)+xy^2L_3(\be)+y^3L_4(\be),
\end{equation}
for certain linear forms $L_j(\be)$.  Similarly we may write
\begin{equation}\label{5c}
F_2(x,y;\be)=x^2Q_1(\be)+xyQ_2(\be)+y^2Q_3(\be),
\end{equation}
where $Q_1,Q_2,Q_3$ are quadratic forms.  We now substitute for
$x^2,xy$ and $y^2$ according to (\ref{inv}), whence
\[F_2(x,y;\be)=Q(x,y)R_1(\be)+Q'(x,y)R_2(\be)+Q''(x,y)R_3(\be),\]
for quadratic forms $R_j(\x)\in\Z_p[\x]$.
Finally, we specify that in this case $\be$ must satisfy the
conditions
\[L_1(\be)=L_2(\be)=L_3(\be)=L_4(\be)=R_2(\be)=R_3(\be)=0.\]

Overall we see that the vector $\be$ must satisfy at most 16 linear conditions
and 8 quadratic conditions.  This is possible when
\[n>V_2(8,16;p)=16+\beta(8;\Q_p).\]
Let us write $r$ for the level of $\be$, and multiply by an
appropriate power of $p$ so that $v(F(\be))=r$.  Clearly the
maximality of $S$ implies that there is at least one vector $\be_i$ of
level $r$.

We begin by examining the possibility that there is just one vector
$\be_i$ of level $r$. Then
\[p^{-r}F(x\be_i+y\be)=ax^4+bxy^3+cy^4\]
for certain $a,b,c\in\Q_p$, by construction.  Moreover we have $a,c\in \Z_p$.
Lemma \ref{li} shows that $\be_i$ and $\be$ are linearly independent,
and then Lemma \ref{mcoe} shows that $b\in\Z_p$, since $F(\x)$ has no
non-trivial zeros.  It follows that $S\cup\{\be\}$ is an
admissible set, contradicting the maximality of $S$.  Hence there
cannot be exactly one vector $\be_i$ of level $r$.

We now suppose that there are two vectors $\be_i,\be_j\in S$ of level
$r$.  The constraints imposed on $\be$ above show that (with
a slight change of notation)
\begin{eqnarray}\label{Quar}
p^{-r}F(x\be_i+y\be_j+z\be)&=&H(x,y)\nonumber\\
&&\hspace{3mm}\mbox{}+DQ(x,y)z^2+(Ex+Fy)z^3+Gz^4,
\end{eqnarray}
where $H(x,y)$ is a binary form with coefficients in $\Z_p$.  Moreover
$\theta(H(x,y))=Ax^4+Bxy^3+Cy^4$.
As usual, Lemma \ref{li} shows that $\be_i,\be_j$ and $\be$ are 
linearly independent. 

We must next prove that $D, E$ and $F$ in (\ref{Quar}) are $p$-adic
integers. We shall argue by contradiction.  Suppose that
\[s:=\min\{v(D),v(E),v(F)\}<0.\]
Then 
\[\theta(p^{-r-s}F(x\be_i+y\be_j+z\be))=
dq(x,y)z^2+(ex+fy)z^3\in\F_p[x,y,z], \]
where at least one of
$d,e$ and $f$ is non-zero.  Here we have recalled that the quadratic
form $Q$ was chosen to be a lift of $q$.  Now unless $e$ and $f$ both
vanish, the point $(0,0,1)$ is a non-singular solution to
$dq(x,y)z^2+(ex+fy)z^3=0$, which therefore lifts to a $p$-adic
solution of $F(x\be_i+y\be_j+z\be)=0$, by Hensel's Lemma.  This
contradicts our assumption that the only $p$-adic zero of $F(\x)$ is
the trivial one.  Hence we must have $e=f=0$ and $d\not=0$.  However
the form $q(x,y)$ was constructed to have distinct linear factors over
$\F_p$, whence $q(x,y)=0$ has a non-singular solution $(a,b)$ say,
leading to a non-singular solution $(a,b,1)$ of $dq(x,y)z^2=0$.  This
again can be lifted to produce a non-trivial solution of $F(\x)=0$.
Thus we have a contradiction unless $D, E$ and $F$ are $p$-adic
integers.

Finally, we conclude that 
\begin{equation}\label{badform}
\theta(p^{-r}F(x\be_i+y\be_j+z\be))
\end{equation}
is of the form
(\ref{gf}) in Lemma \ref{mykey}.  If the form has a non-singular
zero we can apply Hensel's Lemma to produce a
non-trivial solution of $F(\x)=0$.  Thus the only difficulty arises
when the level $r$ is either unsuitable or unacceptable, and 
either $p\equiv 5,\,7\pmod{8}$
with (\ref{badform}) similar to (\ref{bf}), or $p\in\{13,29\}$ with
(\ref{badform}) diagonal.  In the second case computation shows that 
(\ref{badform}) will have a non-singular zero except when it is
similar, after permutation of the variables, to 
$x^4+y^4+2z^4$ (for $p=13$), or $x^4+y^4+z^4$ (for $p=29$).  Of
course, when (\ref{badform}) has a non-singular zero we can produce a
zero of the original form $F(\x)$ via Hensel's Lemma.

We now come to the key step for these remaining cases.  If the level
$r$ is unsuitable we replace $\be_j$ by $\be$ to form a new set $S'$.
Then $S'$ will be admissible, and will have the same size as $S$.
However, since $\theta(p^{-r}F(x\be_i+z\be))$ is similar to
$x^4+2H^2z^4$ when $(\ref{badform})$ is similar to $(\ref{bf})$, we see
that $S'$ has one fewer unsuitable level. This contradicts our
original choice of $S$.

Similarly, if the level $r$ is unacceptable we observe that
\[(x+y)^4+(2x+y)^4+2(x+2y)^4=6x^4+11xy^3+8y^4\]
in $\F_{13}$, and
\[(x+y)^4+(6x+26y)^4+(x+9y)^4=22x^4+10xy^3+2y^4\]
in $\F_{29}$.  Moreover, the form $6x^4+11xy^3+8y^4$ is not similar to
$x^4-4xy^3+3y^4$ over $\F_{13}$, and $22x^4+10xy^3+2y^4$ is not similar to
$x^4-4xy^3+3y^4$ over $\F_{29}$.

Thus, for an unacceptable level, 
there will be a pair of linearly independent
vectors $\be'_i,\be'_j$ in the span of $\{\be_i,\be_j,\be\}$ such that
$p^{-r}F(x\be'_i+y\be'_j)$ has coefficients in $\Z_p$ and such that
$\theta(p^{-r}F(x\be'_i+y\be'_j))$ is of the shape 
$6x^4+11xy^3+8y^4$ or $22x^4+10xy^3+2y^4$ as
appropriate.  For example, if $p=13$ and
\[\theta(p^{-r}F(x\be_i+y\be_j+z\be))=x^4+y^4+2z^4\]
then we set $\be'_i=\be_i+2\be_j+\be$ and $\be'_j=\be_i+\be_j+2\be$.

We now consider the new set $S'$ formed from $S$ by replacing $\be_i$
and $\be_j$ by $\be'_i$ and $\be'_j$.  It is clear that $S'$ will also
be admissible, and that it will have the same size as $S$.  However it
will have one more acceptable level than $S$, and this level will not
be unsuitable.  This again contradicts our original
choice of the set $S$.

This completes the proof of Lemma \ref{Ll}. 

\section{Theorem \ref{odd} for $p=5$}\label{case5}

It remains to consider the case $p=5$.  Here it seems that we cannot
make do by imposing only two quadratic constraints per level, for the
new vector $\be$.  The difficulty revolves around the possibility of a
level $r$ with two vectors $\be_i$ and $\be_j$ for which 
$\theta(p^{-r}F(x\be_i+y\be_j))=2x^4+y^4$ and such that the new vector
$\be$ also has level $r$ and satisfies
$\theta(p^{-r}F(x\be_i+y\be_j+z\be))=2x^4+y^4+z^4$.  We therefore use
an argument in which we impose up to three quadratic constraints for
each level, and this results in the following larger bound for $v_p(4)$.
\begin{lemma}\label{L5}
We have
\[v_4(5)\le 40+\beta(12;\Q_5).\]
\end{lemma}
Thus the case $r=12$ of Lemma \ref{q} gives us the bound for $v_4(5)$ in
Theorem \ref{odd}.

We begin by specifying what we shall mean by an 
``admissible'' set $S$ for $p=5$.  We require the following conditions.
\begin{enumerate}
\item[(i)] $0\le v(F(\be_i))\le 3$ for $1\le i\le m$.
\item[(ii)] For each level $r$ there are at most three vectors
  $\be_i$ of level $r$.
\item[(iii)] The set of all vectors $\be_i$ of a 
given level is linearly independent.
\item[(iv)] If there are exactly two vectors $\be_i$ and $\be_j$ 
of level $r$, with $i<j$,
  then the binary form $p^{-r}F(x\be_i+y\be_j)$ has coefficients in
  $\Z_p$, and 
\[\theta(p^{-r}F(x\be_i+y\be_j))=Ax^4+Bxy^3+Cy^4\]
for certain $A,B,C\in\F_p$
  depending on $i$ and $j$.
\item[(v)] If there are three vectors $\be_i,\be_j$ and $\be_k$ 
of level $r$, with $i<j<k$,
  then the ternary form $p^{-r}F(x\be_i+y\be_j+z\be_k)$ has coefficients in
  $\Z_p$, and $\theta(p^{-r}F(x\be_i+y\be_j+z\be_k))=c(2x^4+y^4+z^4)$
  for some $c\in\F_p$.
\end{enumerate}

When there are exactly two vectors $\be_i$ and $\be_j$ 
of level $r$, with $i<j$, we say that the level is ``suitable'' if
$\theta(p^{-r}F(x\be_i+y\be_j))=c(2x^4+y^4)$ for some $c\in\F_p$, and
otherwise ``unsuitable''.  We choose a set $S$ whose cardinality is
maximal, and having as few unsuitable levels as possible.  As before
we argue by contradiction, assuming that $F(\x)=0$ has only the
trivial solution, and we produce a further
non-zero vector $\be$ satisfying certain constraints, which we now
describe. 

If the set $S$ has no elements of level $r$
there will be no corresponding constraints.  If $S$ has exactly one
element of level $r$ we require one linear and one quadratic
constraint as in the previous cases.  

When $S$ has
exactly two vectors $\be_i$ and $\be_j$ 
of level $r$ we consider the expansions (\ref{5a}), (\ref{5b}) and
(\ref{5c}), and impose the conditions
\[L_1(\be)=L_2(\be)=L_3(\be)=L_4(\be)=Q_1(\be)=Q_2(\be)=Q_3(\be)=0.\]
Finally, when there are three vectors $\be_i,\be_j$ and $\be_k$ 
of level $r$, we write
\begin{eqnarray*}F(x\be_i+y\be_j+z\be_k+w\be)&=&F(x\be_i+y\be_j+z\be_k)
+F_3(x,y,z;\be)w\\
&&\hspace{3mm}\mbox{}+F_2(x,y,z;\be)w^2+F_1(x,y,z;\be)w^3+F(\be)w^4,
\end{eqnarray*}
where each $F_i(x,y,z;\be)$ is bi-homogeneous, of degree $i$ in $(x,y,z)$
and of degree $4-i$ in $\be$.  In particular we have
\[F_3(x,y,z;\be)=\sum_{d+e+f=3}x^dy^ez^fL_{d,e,f}(\be)\]
and
\[F_2(x,y,z;\be)=\sum_{d+e+f=2}x^dy^ez^fQ_{d,e,f}(\be)\]
for certain linear forms $L_{d,e,f}(\be)$ and quadratic forms $
Q_{d,e,f}(\be)$.  In this case we impose 10 linear constraints
\[L_{d,e,f}(\be)=0,\;\;\; \mbox{for all}\;\;\; d,e,f\ge 0\;\;\;
\mbox{with}\;\;\; d+e+f=3,\]
and three quadratic constraints 
\[Q_{2,0,0}(\be)=Q_{0,2,0}(\be)=Q_{0,0,2}(\be)=0.\]

Overall we see that the vector $\be$ must satisfy at most 40 linear conditions
and 12 quadratic conditions.  This is possible when
\[n>V_2(12,40;5)=40+\beta(12;\Q_5).\] 
We suppose that $\be$ has level $r$ and indeed that $v(F(\be))=r$.  As
in \S \ref{large}, if $S$ contains at most one vector $\be_i$ of
level $r$ we get a contradiction, since $S\cup\{\be\}$ will also be
admissible. 

We now consider the possibility that $S$ contains exactly two vectors
$\be_i$ and $\be_j$ of level $r$.  By construction we have
\[p^{-r}F(x\be_i+y\be_j+z\be)=h(x,y)+(dx+ey)z^3+fz^4,\]
where $h$ has coefficients in $\Z_p$ and $f$ is a $p$-adic unit.
Then $\{\be_i,\be_j,\be\}$ must be linearly independent, by Lemma
\ref{li}.  Moreover, by taking $y=0$ and applying Lemma \ref{mcoe}, we
see that $d$ must be in $\Z_p$, and similarly for $e$.  We now apply
the following modification of Lemma \ref{mykey}.
\begin{lemma}\label{mykey51}
Let 
\[h(x,y)=Ax^4+Bxy^3+Cy^4\in\F_5[x,y]\]
be a binary quartic form with $AC\not=0$.
Then for any $D,E,F\in\F_5$ with $F\not=0$, either the form 
\[g(x,y,z):=h(x,y)+Dxz^3+Eyz^3+Fz^4\]
has at least one non-singular zero over $\F_p$, or we can permute the
variables $x,y,z$ to give
\[g(x,y,z)=c(x^4+y^4+z^4)\;\;\mbox{or}\;\;c(2x^4+y^4+z^4)\;\;\mbox{or}\;\;
c(x^4+y^4+dxz^3+3z^4)\]
for certain $c,d\in\F_5-\{0\}$.
\end{lemma}
This may be established by a direct computer check.

If $g(x,y,z)=\theta(p^{-r}F(x\be_i+y\be_j+z\be))$ were to have a 
non-singular zero it could be lifted to a non-trivial
zero of $F(x\be_i+y\be_j+z\be)$ over $\Q_5$, thereby giving a
contradiction.  On the other hand if
\[g(x,y,z)=\theta(p^{-r}F(x\be_i+y\be_j+z\be))=c(x^4+y^4+z^4),\]
then 
$\theta(p^{-r}F(x\be_i+y\be_j))=c(x^4+y^4)$, whence
the level $r$ must have been unsuitable.  In this case we observe that
$g(x,y,x)=c(2x^4+y^4)$.  Thus if we set
$\be'_i=\be_i+\be$ and replace $\be_i$ by $\be'_i$ in $S$, we will produce
a new set $S'$ with
$\theta(p^{-r}F(x\be'_i+y\be_j))=c(2x^4+y^4)$.  It follows that $S'$
has one fewer unsuitable level than $S$, which contradicts our choice
of $S$.  We argue similarly if $g(x,y,z)=
\theta(p^{-r}F(x\be_i+y\be_j+z\be))=c(x^4+y^4+dxz^3+3z^4)$, using the
fact that $g(x,y,2dx)=c(2x^4+y^4)$.  Again we will produce an
admissible set $S'$ with one fewer unsuitable level than $S$, 
contradicting our choice of $S$. Finally, if
$\theta(p^{-r}F(x\be_i+y\be_j+z\be))=c(2x^4+y^4+z^4)$, then we can
take $S'=S\cup \{\be\}$, which will be admissible, since condition (v)
is now satisfied in our definition.  This contradiction shows that $S$
cannot have exactly two vectors of level $r$.

To complete our treatment of the case $p=5$ we examine the situation
in which $S$ has three vectors $\be_i,\be_j$ and $\be_k$ with the same
level $r$ as $\be$. By construction we now have
\[p^{-r}F(x\be_i+y\be_j+z\be_k+w\be)=h(x,y,z)+q(x,y,z)w^2+
(dx+ey+fz)w^3+gw^4,\]
where $h$ has coefficients in $\Z_p$ and $g$ is a $p$-adic
unit. Moreover the quadratic form $q(x,y,z)$ takes the shape
\[q(x,y,z)=axy+bxz+cyz\]
with $a,b,c\in\Q_p$.  As before, the set $\{\be_i,\be_j,\be_k,\be\}$ 
must be linearly independent, by Lemma \ref{li}.  Moreover, by taking
two of $x,y$ and $z$ to vanish, and applying Lemma \ref{mcoe}, we
see that each of $d,e$ and $f$ must be in $\Z_p$.  We proceed to show
that $a,b$ and $c$ are also in $\Z_p$.  Suppose to the contrary that
$v(a)=s<0$, say, with $a=p^sa'$.  Then on setting $z=0$ we have
\[\theta(p^{-r-s}F(x\be_i+y\be_j+w\be))=\theta(a')xyw^2,\]
which has a non-singular zero at $(x,y,w)=(0,1,1)$.  By Hensel's Lemma
we may then derive a nontrivial zero of $F(x\be_i+y\be_j+w\be)$,
contradicting our basic assumption.  Thus $a$ must be a $p$-adic
integer, and similarly for $b$ and $c$.  Finally we apply the
following lemma.
\begin{lemma}\label{mykey52}
Let 
\[H(x,y,z)=2x^4+y^4+z^4\in\F_5[x,y,z].\]
Then for any $A,B,C,D,E,F,G\in\F_p$ with $G\not=0$ the form 
\[g(x,y,z,w):=H(x,y,z)+(Axy+Bxz+Cyz)w^2+(Dx+Ey+Fz)w^3+Gw^4\]
has at least one non-singular zero over $\F_5$.
\end{lemma}
Again this is the result of a computer check.  Lemma \ref{mykey52}
now shows that 
\[\theta(p^{-r}F(x\be_i+y\be_j+z\be_k+w\be))\]
has a
non-singular zero, whence $F(x\be_i+y\be_j+z\be_k+w\be)$ has a
non-trivial zero in $\Q_5$.  Thus $F(\x)$ has a non-trivial zero.
This contradiction establishes Lemma \ref{L5}.

\section{The Proof of Theorem \ref{even}}

The methods employed to prove Theorem \ref{odd} are based on the 
application of Hensel's Lemma to lift zeros of forms defined over
$\F_p$. We have no way to guarantee that the forms we construct will
not be diagonal, in which case there will be no
non-singular zeros over $\F_2$.  Thus it would appear that
the approach is completely inapplicable for $p=2$. Our treatment of
Theorem \ref{even} will therefore be based largely on Wooley's version
of the quasi-diagonalization method.  However we will make extensive 
use of the idea introduced in \S \ref{large}, where we used the newly
constructed $\be$ to alter one of the vectors in $S$, rather than
merely adding $\be$ to $S$.  

Our primary goal in this section is to prove the following bound.
\begin{lemma}\label{evengoal}
We have
\[v_4(2)\le V_3(5,21,56;2).\]
\end{lemma}
The estimate given in Theorem \ref{even} is then an immediate
consequence of (\ref{med5}) in conjunction with Lemmas \ref{q} and 
\ref{phipsi}.

We assume throughout this section that the form
\[F(x_1,\ldots,x_n)\in\Q_2[x_1,\ldots,x_n]\]
is fixed, and that $F$ has
only the trivial 2-adic zero. Given a set 
$S=\{\be_1,\ldots,\be_k\}$ of non-zero vectors in $\Q_2^n$ we shall say that
a non-zero vector $\be\in\Q_2^n$ is ``orthogonal'' to $S$ if
\[F(x_1\be_1+\ldots +x_k\be_k+x\be)=F(x_1\be_1+\ldots +x_k\be_k)
+F(\be)x^4.\]
Thus, by Lemma \ref{li}, if $S$ is linearly independent, then so is
$S\cup\{\be\}$.  The following result tells us when such an $\be$
exists.
\begin{lemma}\label{orth}
If $\#S=k$ and
\[n>V_3(k,\frac{k(k+1)}{2},\frac{k(k+1)(k+2)}{6};2)\]
there is a vector $\be$ orthogonal to $S$.
\end{lemma}
For the proof we observe that we can write
\begin{eqnarray*}
F(x_1\be_1+\ldots +x_k\be_k+x\be)&=&F(x_1\be_1+\ldots +x_k\be_k)
+\sum_{\sum d_i=3}\x^{\lb{d}}F^{(1)}_{\lb{d}}(\be)x\\
&&\mbox{}+\sum_{\sum d_i=2}\x^{\lb{d}}F^{(2)}_{\lb{d}}(\be)x^2
+\sum_{\sum d_i=1}\x^{\lb{d}}F^{(3)}_{\lb{d}}(\be)x^3\\
&&\mbox{}+F(\be)x^4,
\end{eqnarray*}
where the forms $F_{\lb{d}}^{(m)}(\be)$ all have degree $m$ in $\be$.
Thus, in order to ensure that $\be$ is orthogonal to $S$ it suffices
that all the forms $F_{\lb{d}}^{(m)}(\be)$ should vanish for $1\le
m\le 3$ and $\sum_{i=1}^k d_i=m$.  Thus $\be$ must be a simultaneous
zero of a system of $k(k+1)(k+2)/6$ linear forms, $k(k+1)/2$ quadratic
forms, and $k$ cubic forms. The result then follows.

We may construct diagonal forms 
$F(x_1\be_1+\ldots +x_k\be_k)$ by using Lemma \ref{orth} iteratively.  
We then say that the vectors
$\be_1,\ldots,\be_k$ are ``mutually orthogonal''. A convenient
criterion for when such a diagonal form has a non-trivial 2-adic zero
is given by the next lemma. Here we use the notion of the ``level'' of
a vector, as introduced in \S \ref{secC}
\begin{lemma}\label{25}
Let $\be_1,\ldots,\be_5$ be mutually orthogonal, and suppose that
there is at least one vector of each level $r$, for $0\le r\le
3$. Then $F(x_1\be_1+\ldots +x_5\be_5)$ has a non-trivial $2$-adic
zero.

If $\be_1,\ldots,\be_4$ are mutually orthogonal, with exactly
one vector of each level $r\in\{0,1,2,3\}$, 
then $F(\x)$ has a non-trivial $2$-adic
zero providing that
\[n>V_3(4,10,20;2).\]
\end{lemma}
The second statement is an immediate deduction from the first, since
Lemma~\ref{orth} enables us to find a fifth vector $\be_5$ orthogonal to
$\be_1,\ldots,\be_4$.

To prove the first statement we consider diagonal 2-adic forms
$\sum_1^5 c_ix_i^4$. We multiply the form by an appropriate power of
2, re-order the indices, and re-scale the variables by powers of 2, so
that 
\[v(c_1)=v(c_2)=0,\; v(c_3)=1,\; v(c_4)=2,\;\;\mbox{and}\;\; v(c_5)=3.\]
Indeed, dividing the form by $c_1$, we may assume that $c_1=1$.  Since
2 divides $1+c_2$, and $v(c_3)=1$, we can choose $x_3\in\{0,1\}$ so that 
$4|1+c_2+c_3x_3^4$.  By the same reasoning we can then select $x_4\in\{0,1\}$
so that $8|1+c_2+c_3x_3^4+c_4x_4^4$, and $x_5\in\{0,1\}$
so that $16|1+c_2+c_3x_3^4+c_4x_4^4+c_5x_5^4$.  We now set
$x_2=1$ and $A=-\sum_2^5c_ix_i^4$, whence $A\equiv 1\pmod{16}$.  Then
$A$ is a fourth power in $\Z_2$, equal to $x_1^4$, say.  It follows
that $\sum_1^5 c_ix_i^4=0$ with the $x_i$ not all zero, as required.

We now assume that 
\[n>V_3(5,15,35;2),\]
whence successive applications of Lemma \ref{orth} allow us to
construct a mutually orthogonal set $\be_1,\ldots,\be_6$. It follows
from Lemma \ref{25} that not all four levels can be attained by these
vectors, since we are supposing that $F(\x)$ has only the trivial
zero.  We proceed to investigate just what one can say about the
levels of vectors in such a mutually orthogonal set.  The basic
principle we shall use is embodied in the following result.
\begin{lemma}\label{shift}
Let $F(x_1,\ldots,x_n)\in\Q_2[x_1,\ldots,x_n]$ have no non-trivial
$2$-adic zero, and suppose that
\[n>V_3(5,15,35;2).\]
Suppose that the set $\be_1,\ldots,\be_6$ is mutually orthogonal and
that $\be_1,\be_2$ and $\be_3$ all have the same level $r$.  Then
there is a mutually orthogonal set $\be'_1,\be'_2,\be'_3,\be_4,\be_5,\be_6$
in which $\be'_1$ has level $r$ and $\be'_2$ has level $r+1$ (or level
$0$, in case $r=3$).
\end{lemma}

For the proof we assume for simplicity that $r=0$, the other cases
being similar.  Under this assumption we have
\[F(x_1\be_1+\ldots +x_6\be_6)=c_1x_1^4+\ldots +c_6x_6^4\]
with $c_1,c_2$ and $c_3$ being 2-adic units.  It follows that $c_i\equiv\pm
1\pmod{4}$ for $1\le i\le 3$, whence there are two 
indices $1\le i<j\le 3$ such that
$c_i\equiv c_j\pmod{4}$.  In particular we will have $c_i+c_j\equiv
2\pmod{4}$. If $k$ is the third index in $\{1,2,3\}$ we set
$\be'_1=\be_k$ and $\be'_2=\be_i+\be_j$.  Hence
\[F(x_1\be'_1+x_2\be'_2+x_4\be_4+x_5\be_5+x_6\be_6)=
c_kx_1^4+(c_i+c_j)x_2^4+c_4x_4^4+c_5x_5^4+c_6x_6^4,\]
so that $\be'_1$ has level $0$ and $\be'_2$ has level 1.  We complete
the proof by applying Lemma \ref{orth} to obtain an additional
orthogonal vector $\be'_3$.

We may use Lemma \ref{shift} to produce an orthogonal set with a
convenient collection of levels.
\begin{lemma}\label{123}
Let $F(x_1,\ldots,x_n)\in\Q_2[x_1,\ldots,x_n]$ have no non-trivial
2-adic zero, and suppose that
\[n>V_3(5,15,35;2).\]
Then, for an appropriate integer $k$, the form $2^kF(\x)$ has an
orthogonal set $S=\{\be_1,\ldots,\be_6\}$ in which $\be_1$ and $\be_2$
have level $0$, $\be_3$ and $\be_4$ have level $1$, and $\be_5$ and $\be_6$
have level $2$. 
\end{lemma}

We begin the proof by showing that there is an orthogonal set with at
least 3 different levels. 
Lemma \ref{shift} shows that if the vectors in $S$ all 
have the same level then we may replace them by a new set in which at least
two different levels appear.  Suppose now that we have a set
$S$ containing precisely two different levels.  
We may multiply $F$ by a suitable power of 2 so that the two
levels present in our original set $S$ are either 0 and 1 or 0 and 2.  
It is easy to dispose of the latter case, since at least one of the
levels 0 or 2 must occur for three or more vectors $\be_i$.  Suppose
for example that $\be_1,\be_2$ and $\be_3$ have level 0 and that
$\be_4$ has level 2.  Then an
application of Lemma \ref{shift} will produce a set $S'$ containing
vectors $\be'_1$ of level 0, $\be'_2$ of level 1 and $\be_4$ of level 2.

To deal with sets $S$ which have levels 0 and 1 and no others, we
consider such a set $S$ in which the number of vectors of level 1 is
maximal.  If this set has 3 or more vectors of level 0 we may apply 
Lemma \ref{shift} to produce a new set $S'$ with an additional element
of level 1, and this would contradict our assumption unless $S'$ has
3 different levels. On the other hand, if $S$ has 1 or 2 elements of
level 0 then there are 4 or 5 elements of level 1.  Thus we may apply 
Lemma \ref{shift} to produce a set $S'$ with at least one element of level $r$
for $r=0,1$ and 2. Hence we may always obtain an orthogonal set with
at least three different levels.  Of course if there are 4 different
levels then the second assertion of Lemma \ref{25} gives a contradiction.

We now show that if we have an orthogonal set $S$ with 3 different
levels we can derive a new set $S'$ with precisely the levels
specified in Lemma \ref{123}. By appropriate choice of $k$ we may
assume that $S$ has elements of levels 0, 1 and 2.  If the numbers of
elements of these levels are $a,b$ and $c$ respectively we will 
assign a ``score'' $b+3c$ to the set $S$.  We now consider such a
set with the maximum score possible.  If $a\ge 3$ we can apply Lemma
\ref{shift} to $S$ to obtain a set $S'$ with score $b'+3c'$, and with
$b'\ge b+1$ and $c'\ge c$.  Thus $S'$ would have a larger score than $S$.
Similarly if $b\ge 3$ we can apply Lemma \ref{shift} to produce a set
$S'$ with $b'\ge b-2$ and $c'\ge c+1$.  Again this shows that $S'$
would have a larger score than $S$.  Finally, if $S$ has elements of
levels 0, 1 and 2, and has $c\ge 3$, Lemma \ref{shift} will produce a
set $S'$ containing all four levels.  However this is impossible since
the second part of Lemma \ref{25} would then show that $F(\x)$ has a
non-trivial zero.  Thus our set $S$ can only have $a=b=c=2$, as required.

Before completing the proof of Lemma \ref{evengoal} we observe that
one can investigate orthogonal sets of size 7 in much the same way as
we have done here for sets of size 6.
In this case repeated use of Lemma \ref{shift} will always eventually
lead to an orthogonal set containing vectors of all four levels, so
that Lemma \ref{25} can be applied.  Hence we will have
\[v_4(2)\le V_3(6,21,56;2).\]
However Lemma \ref{evengoal} improves on this somewhat.

To establish Lemma \ref{evengoal} we start from the set $S$ constructed in
Lemma \ref{123}, so that
\[2^kF(x_1\be_1+\ldots +x_6\be_6)=c_1x_1^4+\ldots+c_6x_6^4\]
with $v(c_1)=v(c_2)=0$, $v(c_3)=v(c_4)=1$ and $v(c_5)=v(c_6)=2$.  We
proceed to find a further vector $\be$ which is ``nearly'' orthogonal
to $\{\be_1,\ldots,\be_6\}$.  Specifically we shall require that
\[2^kF(x_1\be_1+\ldots +x_6\be_6+x\be)=2^kF(x_1\be_1+\ldots
+x_6\be_6)+Ax_1x^3+Bx^4\]
for some $A,B\in\Q_2$.
An argument completely analogous to that used for Lemma~\ref{orth}
shows that this is possible with $\be\not=\b{0}$, providing that we
can satisfy simultaneously 56 linear constraints,
21 quadratic constraints and 5 cubic constraints. Hence
$n>V_3(5,21,56;2)$ suffices.  By Lemma \ref{li} the set
$\be_1,\ldots,\be_6,\be$ will be linearly independent.  Moreover,
since we are assuming that $F(\x)$ has no non-trivial zero, we will
have $B\not=0$.  Thus, by re-scaling $\be$ by a power of 2, we may
assume that $v(B)=0,1,2$ or 3.

We now observe that for any $a\in\Q_2$ the set
$S_a=\{a\be_1+\be,\be_2,\be_3,\be_4,\be_5,\be_6\}$ will be orthogonal, and
certainly contains vectors of levels 0,1 and 2. Suppose that
$a\be_1+\be$ has level $\lambda$. We cannot have $\lambda=3$, since
then Lemma \ref{25} would produce a non-trivial zero of $F(\x)$.  If
$\lambda=2$ then $S_a$ has 1 element of level 0; it has 2 elements of level
1; and 3 elements of level 2.  In this case an application of Lemma
\ref{shift} will produce a new orthogonal set $S_a'$ containing
elements of all four levels, which is impossible by Lemma \ref{25}.
Similarly if $\lambda=1$ then $S_a$ has 1 element of level 0; there are 3
elements of level 1; and 2 elements of level 2.  This time Lemma
\ref{shift} yields a set $S_a'$ with at least one element of each of
the levels 0 and 1, and at least 3 elements of level 2.  Thus a second
application of the lemma gives us a set $S_a''$ containing all four
levels, which again gives a contradiction via Lemma \ref{25}.

There remains the possibility that $\lambda=0$ for every choice of $a$.
In particular, taking $a=0$, we see that $B$ must be a 2-adic unit.
Lemma \ref{mcoe} then shows that $A\in\Z_2$.  We now consider
the polynomial
\[f(x)=2^kF(x\be_1+\be_2+\be)=c_1x^4+Ax+B+c_2.\]
If $A$ is a 2-adic unit then $\theta(f(x))=x^4+x$ which has a
non-singular zero in $\F_2$, at $x=1$.  By Hensel's Lemma this would
produce a zero of $f(x)$ in $\Z_2$, and hence a non-trivial zero of
$F(\x)$.  We therefore conclude that $2|A$.  Thus $F(a\be_1+\be)$ must
be even whenever $a$ is a 2-adic unit, and since $a\be_1+\be$ has
level zero we deduce that $16|F(a\be_1+\be)$. Taking $a=\pm 1$ we find
that 
\[c_1\pm A+B\equiv 0\pmod{16},\]
so that $8|A$.  We now choose
$t=0$ or 2 such that $32|c_1+B+A+c_2t^4$, and consider the polynomial 
\[g(x)=2^kF(x\be_1+t\be_2+\be)=c_1x^4+Ax+B+c_2t^4.\]
By construction we have $2^5|g(1)$ and
\[g'(1)=4c_1+A\equiv 4\pmod{8},\]
so that $2^3\nmid g'(1)$.  It follows from Hensel's Lemma that $g(x)$
has a zero in $\Z_2$, and hence that $F(\x)$ has a non-trivial zero in
$\Q_2$.  This completes the proof of Lemma \ref{evengoal}.

\bigskip
\bigskip

Mathematical Institute,

24--29, St. Giles',

Oxford

OX1 3LB

UK
\bigskip

{\tt rhb@maths.ox.ac.uk}

\end{document}